\documentstyle[12pt,righttag]{amsart}
\numberwithin{equation}{subsection}
\theoremstyle{plain}
  \newtheorem{th}{Theorem}[subsection]
  \newtheorem{pr}[th]{Proposition}
  \newtheorem{lem}[th]{Lemma}
  \newtheorem{cor}[th]{Corollary}
\theoremstyle{definition}
  \newtheorem{df}[th]{Definition}
\theoremstyle{remark}
 \newtheorem{rem}[th]{Remark}

\def\e{{\epsilon}}
\def\g{{\frak{sl}}}
\def\h{{\frak h}}

\def\l{{\ell}}

\def\n{{\frak n}}

\def\t{{{\frak t}}}

\def\ob{{obj\,}}
\def\proj{{p}}
\def\act{{\pi}}

\def\scsc{\scriptscriptstyle}
\def\rkH{{\ell}}
\def\rkg{{n}}

\def\Vn{{{V^{\*\rkH}_\rkg}}}
\def\Vnn{{{V^{\*\rkH}_\rkH}}}
\def\Para{{\cal S}}

\def\C{{\Bbb C}}

\def\K{{\cal M}}
\def\Sym{{\frak S}}
\def\L{{\cal L}}

\def\O{{\cal O}}

\def\W{{W}}

\def\Z{{\Bbb Z}}
\def\+{\mathop{\oplus}}
\def\*{\mathop{\otimes}} 
\def\ev{{\text{\rom{ev}}}}

\def\half{{\frac{1}{2}}}
\def\pair{{\lm,\mu}}
\def\one{{\text{\rom{\bf 1}}}}

\def\g{{{\frak{sl}}}}

\def\lm{\lambda}
\def\al{\alpha}
\def\Rep{{{\cal R}}}
\def\ch{^\vee}

\def\bra{{\langle}}
\def\ket{{\rangle}}
\def\End{{\text{\rom{End}}}\,}
\def\Hom{{\text{\rom{Hom}}}\,}
\def\dim{{\text{\rom{dim}}}\,}

\def\int{{{\Z}}}
\def\Ker{{\text{\rom{Ker}}}}

\begin{document}
\title[degenerate affine hecke algebras]{
Lie algebras and degenerate affine Hecke algebras of type $A$} 
\renewcommand{\thefootnote}{\fnsymbol{footnote}}
\maketitle
\begin{center}
Tomoyuki~Arakawa

\medskip
{\it Graduate School of Mathematics,
Nagoya University, Japan}

\medskip
and
\medskip

Takeshi~Suzuki\footnote[2]{
Supported by JSPS the Research Fellowships 
for Young Scientists.}

\medskip
{\it Research 
Institute for Mathematical Sciences, Kyoto University, Japan}
\end{center}
\begin{abstract}
We construct a family of exact functors from the 
Bernstein-Gelfand-Gelfand category $\O$
of $\g_\rkg$-modules to 
the category of finite-dimensional representations of 
the degenerate affine Hecke algebra $H_\rkH$ of $GL_\rkH$.
These functors transform Verma modules to standard modules or zero, 
and simple modules to simple modules or zero.
Any simple $H_\rkH$-module can be thus obtained.
\end{abstract}
\begin{center}
{\bf Introduction}
\end{center}

\medskip
The classical Frobenius-Schur-Weil duality 
gives a remarkable correspondence between
the category of finite-dimensional representations of
the symmetric group $\Sym_\rkH$ and
the category of finite-dimensional representations of
the special (or general) linear group $SL_\rkg$.
Its generalizations have been 
studied in e.g.
\cite{CP,Ch;new interpretation,Dr,Ji,VV}
where $\Sym_\rkH$ is replaced
by other algebras, e.g. the Hecke algebras, the 
(degenerate) affine Hecke algebras
or the double affine Hecke algebras,
and $SL_n$ is replaced by the corresponding quantum groups.

In this paper, we present a new direction in
generalizing the classical duality.
Let $\O(\g_n)$ denote
the BGG category of representations
of the complex Lie algebra $\g_\rkg$,
and let
$\Rep(H_\rkH)$ denote the category
of finite-dimensional representations of
the degenerate (or graded) 
affine Hecke algebra $H_\rkH$ of $GL_\rkH$.
To each weight $\lm$ of $\g_\rkg$
such that $\lm+\rho$ is dominant integral
(where $\rho$ is the half sum of the positive roots),
we associate
a functor $F_\lm$ from $\O(\g_\rkg)$ to $\Rep(H_\l)$. 
When we take $\lm=0$ and restrict the functor $F_0$ to the 
category of 
finite-dimensional representations of $\g_\rkg$,
we obtain the classical duality.

To be more precise,
let $V_\rkg=\C^\rkg$ be the vector representation of $\g_\rkg$
and $M(\lm)$ the highest weight Verma module with 
highest weight $\lm$.
For $X\in \ob\O(\g_\rkg)$,
we construct in \S \ref{ss;action} an action of  $H_\rkH$
on $X\*\Vn$ commuting with the $\g_\rkg$-action.
This induces an $H_\rkH$-action
on a certain subquotient 
$F_\lm(X)$ (see \eqref{eq;definition}) of $X\*\Vn$.
When  $\lm+\rho$ is dominant integral,
$F_\lm(X)$ is identified with
$\Hom_{\g_\rkg}(M(\lm),X\*\Vn)$.
(This space is an analogue of the space of the conformal blocks
in the conformal field theory. 
See \cite{AST}.)

Under the assumption that $\lm+\rho$ is dominant 
integral,
we prove that

\smallskip
(1) $F_\lm$ is exact,

\smallskip
(2) $F_\lm$ sends a Verma module to a ``standard module'' 
unless it is zero.

\smallskip
\noindent
Here the standard module is an induced module 
from a certain one dimensional representation of 
a parabolic subalgebra of $H_\rkH$, 
and it has
a unique simple quotient.

Moreover, in the case $\rkg=\rkH$,
we prove that

\smallskip
(3) $F_\lm$ sends a simple module to a simple module unless 
it is zero.

\smallskip
\noindent
We should remark that
our proof of (3) relies on the
formula \eqref{eq;two_multiplicities} in Appendix.
This formula is a consequence of the 
fact that the irreducible decompositions
of the standard modules of $H_\rkH$
(\cite{Gi})
and those of the Verma modules 
(\cite{BB,BK}) are both described by the
Kazhdan-Lusztig polynomials.

We also determine
when the image of the functor is non-zero
(Theorem \ref{th;Verma_to_standard}, 
Theorem \ref{th;simple_to_simple}).
Furthermore, it follows from Zelevinsky's results in
\cite{Zel;induced}
that
any simple $H_\rkH$-module with ``integral weights'' 
is of the form $F_\lm(L)$
for some weight $\lm$ such that $\lm+\rho$ is dominant integral
and some simple $\g_\rkH$-module $L$ 
(see Corollary \ref{cor;anysimple} for the precise statement).

\medskip
\noindent{\bf Acknowledgment.}
Our work is motivated by the joint research
with A.~Tsuchiya  presented in \cite{AST}.
We are grateful to him
for discussions and his encouragement.
We thank T.~Tanisaki
for drawing
our attentions to \cite{Zel;p-adic KL,Zel;two remarks}. 
Thanks are also due to I.~Cherednik and M.~Kashiwara
for valuable comments.
\section{Basic Definitions}
\subsection{Root data}
Let $\t_{\rkg}$ be an $\rkg$-dimensional complex vector 
space with the
basis $\{\e\ch_i\mid i=1,\dots,\rkg\}$ and 
the inner product defined by $(\e\ch_i|\e\ch_j)=\delta_{ij}$.
Let $\t^*_{\rkg}$ be its dual and $\{\e_i\}$ be the dual basis of 
$\{\e\ch_i\}$,
which are orthonormal basis with respect to 
the induced inner product: $(\e_i|\e_j)=\delta_{ij}$. 
The natural pairing between $\t_{\rkg}$ and $\t^*_{\rkg}$ 
will be denoted by
$\bra\, ,\,\ket:\t^*_{\rkg}\times \t_{\rkg}\to \C$.
Put 
$\h^*_{\rkg}=\{\sum_{i=1}^{\rkg} \lm_i\e_i\in 
\t_{\rkg}^*\mid \sum_{i=1}^{\rkg}\lm_i=0\}$
and $\h_{\rkg}=\t_{\rkg}/\C\e\ch$, 
where $\e\ch=\sum_{i=1}^{\rkg}\e\ch_i$, 
so that $\h_{\rkg}$ and $\h^*_{\rkg}$ are dual to  each other.
Define the roots and the simple roots by
$\al_{ij}=\e_i-\e_j$ $(i\neq j)$, $\al_i=\al_{i\, i+1}$ 
respectively and
put 
\begin{align}
  R_{\rkg}&=\{\al_{ij}\mid  1\leq i\neq j\leq \rkg\},\\
  R_{\rkg}^+&=\{\al_{ij}\mid  1\leq i<j\leq \rkg\},\\
  \Pi_{\rkg}&=\{\al_i\mid i=1,\dots,\rkg-1\},
\end{align}
then $R_{\rkg}\subseteq \h^*_{\rkg}$ is a root system of type 
$A_{\rkg-1}$.
Define the coroots and the simple coroots
by $\al\ch_{ij}=\e\ch_i-\e\ch_j$ $(i\neq j)$ and $h_i=\al\ch_{ii+1}$,
respectively.

Let $W_{\rkg}\subset \rom{GL}(\t^*_{\rkg})$ be the Weyl group 
associated to
the above data, which is by definition generated by
the reflections $s_{\al}$ $(\al\in R_{\rkg})$ 
defined by 
\begin{equation}
  s_\al(\lm) =\lm-  \bra\lm, \al\ch\ket \al \quad (\lm\in\t^*_{\rkg}).
\end{equation}
We often write $s_{\al_{ij}}=s_{ij}$.
Observe that $W_{\rkg}$
 preserves $\h^*_{\rkg}\subseteq\t^*_{\rkg}$
and $W_{\rkg}$ is isomorphic to the symmetric group $\Sym_{\rkg}$.
We often use another action of $W_\rkg$ on $\h^*_\rkg$,
which is given by
\begin{equation}
  \label{eq;dot_action}
  w\circ\lm=w(\lm+\rho)-\rho\quad (w\in W_\rkg,\ \lm\in\h_\rkg^*),
\end{equation}
where $\rho=\frac{1}{2}\sum_{\al\in R^+_\rkg}\al$.
Put
\begin{align}
  Q_{\rkg}&=\+_{i=1}^{\rkg-1}\Z\al_i\subseteq \h_{\rkg}^*,\label{eq;Q}
  \quad Q_{\rkg}^+=\+_{i=1}^{\rkg-1}\Z_{\geq 0}\al_i, \label{eq;Q_+}\\
  P_{\rkg}&=\{\lm\in\h_{\rkg}^*\mid \bra \lm,h_i\ket\in \Z \quad 
  \hbox{for all } i=1,\dots,\rkg-1\},\label{eq;P}\\
  P_{\rkg}^+&=\{ \lm\in \h^*_\rkg\mid  \bra \lm,h_i\ket\in 
\Z_{\geq 0} \quad 
  \hbox{for all } i=1,\dots,\rkg-1\}.\label{eq;P^+}
\end{align}
An element of $P_\rkg$ (resp. $P^+_\rkg$)
is called an integral (resp. dominant integral)
weight.
\subsection{Lie algebras of type $A$}
Let $\g_{\rkg}$ be the complex 
Lie algebra of type $A_{\rkg-1}$.
We introduce an invariant inner product on $\g_\rkg$ by 
 $(x| y)=\text{\rom{tr}}_{\C^\rkg}(xy)$,
where $\C^\rkg$ denotes the vector representation of $\g_\rkg$.
We identify
$\h_{\rkg}$ (introduced in the previous subsection)
and the Cartan subalgebra of $\g_{\rkg}$
as inner product spaces.
Let $\g_{\rkg}=\n_+\+\h_{\rkg}\+\n_-$ be the triangular 
decomposition with
$\n_\pm=\+_{\al\in R_{\rkg}^+}(\g_\rkg)_{\pm\al}$, 
where $(\g_\rkg)_\al$ 
denotes
the root space corresponding to $\al\in R_\rkg$.
We choose a set of root vectors $\{e_\al\in (\g_\rkg)_\al\mid 
\al\in R_{\rkg}\}$ 
such that $(e_\al | e_{-\al})=1$ holds for all $\al\in R_{\rkg}^+$.
Let $\{h^i\}_{i=1,\dots,\rkg-1}$ be the dual basis 
of the coroots $\{h_i\}_{i=1,\dots,\rkg-1}$ in $\h_\rkg$, so that
$(h^i|h_j)=\delta_{ij}$.
We define the special elements of $\g_{\rkg}\*\g_{\rkg}$ by
\begin{align}
  r&=\half\sum_{i=1}^{\rkg-1} h^i\* h_i+\sum_{\al\in R_{\rkg}^+}
e_\al\* e_{-\al},
\label{eq;r-matrix}\\ 
 \Omega&=\sum_{i=1}^{\rkg-1}h^i\* h_i+\sum_{\al\in R_{\rkg}^+}
\left(e_\al\* e_{-\al}+e_{-\al}\* e_\al\right),\label{eq;Omega}
\end{align}
which will be used later.
For an $\h_{\rkg}$-module $X$ and $\lm\in\h^*_{\rkg}$, put
\begin{align}
 X_\lm         &=\{ v\in X\mid hv=\bra \lm,h\ket v 
\text{ for all }h\in\h_{\rkg}\}\label{eq;weight_space}\\
 P(X)          &=\{\lm\in\h^*_{\rkg}\mid X_\lm\neq 0\}\label{eq;weight}
\end{align}
The space $X_\lm$
is called the weight space of weight $\lm$
and 
an element of $P(X)$ is called a weight of $X$.
\subsection{The BGG category $\O$}
Let $U(\g_\rkg)$ denote the universal enveloping algebra of 
$\g_\rkg$.
Let $\O(\g_\rkg)$ denote the category 
whose objects are those $\g_{\rkg}$-modules
$X$ such that

\medskip
(i) $X$ is finitely generated over $U(\g_{\rkg})$, 

\smallskip
(ii) $X$ is $\n_+$-locally finite i.e. 
$\dim_\C U(\n_+)v <\infty$ for each $v\in X$,

\smallskip
(iii) 
 $X=\+_{\lm\in\h^*_{\rkg}}X_{\lm}$. 

\medskip\noindent
The morphisms of $\O(\g_\rkg)$ 
are, by definition, all the $\g_{\rkg}$-homomorphisms. 
The category $\O(\g_\rkg)$ is closed under the operations such as
taking submodules, forming quotient modules,
finite direct sums, and tensor products with finite-dimensional modules.
For $\lm\in\h^*_{\rkg}$
let $M(\lm)$ denote the highest weight 
Verma module with highest weight $\lm$.
The unique simple quotient of $M(\lm)$ will be denoted by $L(\lm)$.
The modules $M(\lm)$ and $L(\lm)$ are objects of $\O(\g_\rkg)$.

Let $\chi_{\lm}:Z(U(\g_{\rkg}))\to \C$ denote
the infinitesimal character of $M(\lm)$
(i.e. $zv=\chi_\lm(z)v$ for all $z\in Z(U(\g_\rkg))$ and $v\in M(\lm)$). 
It is known that $\chi_\lm= \chi_\mu$
if and only if $\lm= w\circ\mu$ for some $w\in W_\rkg$.
Let $[\lm]$ denote the orbit $W_\rkg\circ\lm$ and
put $Z_\lm=\Ker\chi_\lm\subset Z(U(\g_\rkg))$.
Define the full subcategory $\O(\g_\rkg)_{[\lm]}$ of $\O(\g_\rkg)$
by
$$\ob \O(\g_\rkg)_{[\lm]}=\{X\in 
\ob \O(\g_\rkg)\mid (Z_\lm)^k X=0 \text{ for some } k\}.$$ 
Then any $X\in\ob\O(\g_\rkg)$ admits a decomposition
\begin{equation}
  \label{eq;decomposition}
X=\+_{[\lm]}X^{[\lm]}  
\end{equation}
such that $X^{[\lm]}\in \ob\O(\g_\rkg)_{[\lm]}$.
The correspondence $X\mapsto X^{[\lm]}$
gives an exact functor on $\O(\g_\rkg)$.
\subsection{$\n$-homology of $\g_\rkg$-modules}
We proof some facts on the zero-th $\n_-$-homology space 
\begin{align}
  H_0(\n_-,X)&=X/\n_-X
\end{align}
of 
 a  $\g_\rkg$-module $X$.
The space $ H_0(\n_-,X)$ has a natural $\h_\rkg$-module structure, 
and
it is known that
\begin{align}
P(H_0(\n_-,X))\subseteq W_\rkg\circ\mu,
\end{align}
for  $X\in\O(\g_n)_{[\mu]}$ ($\mu\in\h_\rkg^*$).
Hence we have
\begin{lem}\label{lem;weight}
 For any $X\in\ob\O(\g_n)$ and $\lm\in\h^*_\rkg$, we have
\begin{align*}
H_0(\n_-,X)_{\lm}&=H_0(\n_-,X^{[\lm]})_{\lm},
\end{align*}
where $X^{[\lm]}$ is defined by the decomposition 
\eqref{eq;decomposition}.
\end{lem}
Therefore we have the natural surjection 
\begin{equation}
  \label{eq;maps}
(X^{[\lm]})_\lm\rightarrow
  H_0(\n_-,X)_{\lm}.
\end{equation}
\begin{pr}\label{pr;key_proposition}
Let $\lm+\rho\in P^+_\rkg$. Then the above map \eqref{eq;maps} 
is bijective.
\end{pr}
\begin{pf}
By Lemma \ref{lem;weight}, 
it is enough to show that
$\left(\n_- X \right)_{\lm}=0$
assuming $X\in \ob\O(\g_n)_{[\lm]}$.
To show this, note that 
\begin{equation}
  \label{eq;wt}
  P(X)\subseteq  \{\lm-\beta\mid \beta\in Q^+_{\rkg}\},
\end{equation}
where
we used the assumption $\lm+\rho\in P^+_\rkg$.
Hence we have
 $ P(\n_- X)\subseteq  \{\lm-\beta\mid \beta\in 
Q^+_{\rkg}\backslash\{0\}\},$
which implies $\lm\notin P(\n_- X)$ as required.
\end{pf}
\begin{rem}\label{rem;key_remark}
(i) There also exists a canonical bijection
$$
\Hom_{\g_n}(M(\lm),X)\stackrel{\sim}{\rightarrow} 
(X^{[\lm]})_\lm$$
if $\lm+\rho\in P_n^+$.

\smallskip\noindent
(ii) Similar arguments prove that Proposition 
\ref{pr;key_proposition}
holds
for any $\lm\in\h_n^*$ such that $\bra\lm+\rho,\al\ch\ket
\notin\Z_{< 0}$ for all $\al\in R_n^+$.
But we do not need this.
\end{rem}
\subsection{Degenerate affine Hecke algebras}
Let $\rkH\in \Z_{\geq 1}$.
Let $S(\t_\rkH)$ be the symmetric algebra of $\t_{\rkH}$, 
which is isomorphic to
the polynomial ring over $\t^*_{\rkH}$.
\begin{df}
The {\it degenerate $($graded$)$ affine Hecke algebra} $H_{\rkH}$ 
of $GL_{\rkH}$ is the unital associative algebra over $\C$ defined by
the following properties:

\medskip
(i) As a vector space, $H_{\rkH}\cong \C[W_{\rkH}]\* S(\t_\rkH)$, 
where $\C[W_\rkH]$ denotes
the group algebra of $W_{\rkH}$.

\smallskip
(ii) The subspaces $\C[W_{\rkH}]\*\C$ and $\C\*S(\t_\rkH)$ are
subalgebras of $H_{\rkH}$ (their images 
will be identified with $\C[W_{\rkH}]$ and $S(\t_\rkH)$ respectively).

\smallskip
(iii) The following relations hold in $H_{\rkH}$: 
\begin{equation}
  s_\al\xi-s_\al(\xi)s_\al=-\bra \al,\xi\ket\quad
(\al\in \Pi_{\rkH},\ \xi\in\t_{\rkH}).
\end{equation}
\end{df}
The following two lemmas are well known.
\begin{lem}
  \label{lem;center}
The center $Z(H_{\rkH})$ of $H_\rkH$ is 
$$S(\t_\rkH)^{W_{\rkH}}:=\{p\in S(\t_\rkH)\subset H_\rkH
\mid w(p)=p\text{ for all }w\in W_{\rkH}\}.$$
\end{lem}
\begin{lem}
  \label{lem;evaluation}
There exists a unique algebra homomorphism 
$\ev:H_{\rkH}\to \C[W_{\rkH}]$
$($called the evaluation homomorphism$)$
such that 
$$  \ev(w)=w\ (w\in W_{\rkH}),\ \ 
  \ev(\e\ch_i)=\sum_{1\leq j<i}s_{ji}\ (i=1,\dots, \rkH).$$
\end{lem}
\section{exact functor $F_\lm$}
\subsection{}
Let
$V_{\rkg}:=\C^\rkg$ be the vector representation of $\g_{\rkg}$
and let 
$u_i\in V_\rkg$ 
 $(i=1,\dots,\rkg)$ be the vector with only non-zero entry 1 
in the $i$-th component.
Let $\rkH$ be another positive integer.
For $X\in \ob \O(\g_n)$, we
regard $X\* \Vn$ as 
a $\g_{\rkg}$-module.
For $\lm\in \h^*_{\rkg}$,
we put 
\begin{equation}\label{eq;definition}
  F_\lm(X)      =H_0(\n_-,X\* \Vn)_\lm.
\end{equation}
\begin{pr}\label{pr;exact}
Let $\lm+\rho\in P_n^+$. Then
$F_\lm$ is an exact functor from $\O(\g_n)$
to the category of finite-dimensional vector spaces.
\end{pr}
\begin{pf}
Follows from  Proposition \ref{pr;key_proposition},
because $(\cdot)\*\Vn$, $(\cdot)^{[\lm]}$ and $(\cdot)_{\lm}$ 
are all exact functors.
\end{pf}
\begin{rem}
The space $F_\lm(X)$   for  $\lm\in P_\rkg^{+}$
has been studied by
 Zelevinsky \cite{Zel;resolvant}.
He proved that  $F_\lm$ transforms the BGG resolution of a 
finite-dimensional simple $\g_{\rkg}$-modules
to an exact sequence.
\end{rem}
\subsection{$H_\rkH$-action}\label{ss;action}
We shall define an action of $H_\rkH$ on $F_\lm(X)$ as follows.
For $i=0,1,\dots,\rkH$, define 
$\act_i:\g_{\rkg}\to \g^{\*\rkH+1}_{\rkg}$ by 
$\act_i(g)=1^{\* i}\* g\* 1^{\* \rkH-i}$.
We let $\C[\W_{\rkH}]$ act on $X\*\Vn$ naturally 
via permutations of components of $\Vn$.
The image of $w\in W_{\rkg}$ in $\End_\C(X\* \Vn)$
will be denoted by the same symbol.
Note that as operators on  $X\*\Vn$, the following equality holds:
\begin{equation}
  s_{{ij}}=\Omega_{ij}+\frac{1}{\rkg}\quad (1\leq i \neq j\leq \rkH),
\end{equation}
where $\Omega_{ij}=(\pi_i\*\pi_j)(\Omega)$ and
$\Omega$ is as in \eqref{eq;Omega}.
Consider the following operators on $X\* \Vn$:
\begin{align}\label{eq;y_i}
  y_i=\Omega_{0i}+
\sum_{1\leq j<i}\left(\Omega_{ji}+\frac{1}{\rkg}\right)
+\frac{\rkg-1}{2}
\quad (i=1,\dots,\rkH).
\end{align}
\begin{lem}
As operators on $X\*\Vn$, the following equations holds:
\begin{align}
[y_i,y_j]&=0\quad (i,j=1,\dots, \rkH)\label{eq;commutativity},\\
s_iy_j-y_{s_i(j)}s_i&=-\bra\al_i,\e\ch_j\ket\quad (i=1,\dots, \rkH-1,
j=1,\dots, \rkH)\label{eq;with_simple_reflections}.
\end{align}
\end{lem}
{\it Proof.}
  The first equality \eqref{eq;commutativity} 
follows easily from the following relations:
$$ [\Omega_{ij}+\Omega_{jk},\Omega_{ik}]=0,\quad
[\Omega_{ij},\Omega_{km}]=0,$$
where $i,j,k,m\in\{1,\dots,\rkH\}$ are all distinct.
To show \eqref{eq;with_simple_reflections},
recall the homomorphism $\ev$ in Lemma \ref{lem;evaluation},
with which the operator $y_j$  can be written as
\begin{equation}
  \label{eq;with_evaluation}
  y_j=\Omega_{0j}+\ev(\e\ch_j)+\frac{\rkg-1}{2}.
\end{equation}
Now the equality \eqref{eq;with_simple_reflections}
follows by using
$s_i\Omega_{0j}=\Omega_{0\, s_i(j)}s_i \quad (i=1,\dots,\rkH-1).
\quad\Box$
\begin{th}
\rom{(i)} $($\text{{\rom{c.f.}}} \cite{Ch;monodromy}$)$.
 There exists a unique homomorphism
$$H_\rkH\to \End_{\g_\rkg}(X\*\Vn)$$
 such that
  \begin{align*}
    \e\ch_i&\mapsto y_i\quad (i=1,\dots,\rkH),\\
    w      &\mapsto w\quad (w\in W_{\rkH}).
  \end{align*}
\rom{(ii)} The above homomorphism
induces an action of $H_{\rkH}$ on $F_\lm(X)$.
\end{th}
Evidently the correspondence $X\mapsto F_\lm(X)$ defines
a functor from the category $\O(\g_n)$ to the category
$\Rep(H_\l)$
of finite-dimensional $H_{\rkH}$-modules.
Some remarks about this functor
are in the sequel.
\begin{rem}
(i) 
The above construction of the functor $F_\lm$
arose from \cite{AST}, in which
representations of the {\it degenerate double affine Hecke algebra}
are constructed from representations of the affine Lie
algebra $\widehat{\g}_n$
using the Knizhnik-Zamolodchikov connections 
in
the conformal field theory.

\smallskip
\noindent
(ii) Let us consider the case $\lm=0$.
Then for any  $X\in\ob\O(\g_n)$,
the action of $H_\rkH$ on $F_0(X)$
factors the evaluation homomorphism 
in Lemma \ref{lem;evaluation}.
Namely, $F_0$ is regarded as
a functor from $\O(\g_\rkg)$ to 
the category of finite-dimensional representations of $W_\rkH$.
Restricting the functor $F_0$ to the category of finite-dimensional 
representations of $\g_\rkg$,
we obtain
the classical Frobenius-Schur-Weil duality
(see \cite{Zel;resolvant}).
\end{rem}
\section{Images of Verma modules and simple modules}
\subsection{Standard modules and their simple quotients}
\label{ss;standard}
We will determine explicitly
how Verma modules and their simple quotients in $\O(\g_\rkg)$
are transformed by the functor $F_\lm$.
We first introduce some $H_\rkH$-modules.
A pair $[a,b]$ 
of complex numbers such that $b-a+1\in \Z_{\geq 0}$
is called {\it a segment}. 
For a segment $[a,b]$ such that $b-a+1=\rkH$, there exists
 a unique one-dimensional representation 
$\C_{[a,b]}=\C \one_{[a,b]}$ of $H_{\rkH}$ (we put $H_0=\C$ for
convenience)
such that 
\begin{align}
  \label{eq;condition1}
  w \one_{[a,b]} &=\one_{[a,b]}\quad (w\in W_{\rkH}), \\
  \label{eq;condition2}
  \e\ch_i\one_{[a,b]}&= (a+i-1)\one_{[a,b]}\quad (i=1,\dots,\rkH).
\end{align}
Let $\Delta:=([a_1,b_1],\dots,[a_k,b_k])$
be an ordered sequence of segments 
such that $b_i-a_i+1=\l_i$
and $\rkH=\sum_{i=1}^{k}\rkH_i$.
Regard $H_{\rkH_{1}}\* H_{\rkH_{2}}\*\cdots \* H_{\rkH_{k}}$ as 
a subalgebra of $H_{\rkH}$.
Define an $H_{\rkH}$-module $\K(\Delta)$
by
\begin{equation}
  \label{eq;standard}
  \K(\Delta)=H_{\rkH}\*_{H_{\rkH_{1}}\*\cdots \*H_{\rkH_{k}}}
(\C_{[a_1,b_1]}\*\cdots \*\C_{[a_k,b_k]})
\end{equation}
Evidently $\K(\Delta)$ is a cyclic module with a cyclic weight vector 
\begin{equation}\label{eq;one}
  \one_{\Delta}:=
\one_{[a_1,b_1]}\*\cdots\* \one_{[a_k,b_k]},
\end{equation}
whose weight $\zeta_\Delta$ is given by
\begin{equation}\label{eq;weight_of_one}
  \bra\zeta_\Delta,\e\ch_j\ket=a_i+j-\sum_{k=1}^{i-1}\l_k-1\quad
\text{for}\quad\sum_{k=1}^{i-1}\l_k  < j\leq\sum_{k=1}^{i}\l_k.
\end{equation}
It is also obvious that 
$\K(\Delta)\cong \C[\W_{\rkH}/
(\W_{\rkH_{1}}\times\cdots\times\W_{\rkH_k})]$
as a $\C[\W_{\rkH}]$-module. In particular
\begin{equation}
 \label{eq;dim_of_standard}
 \dim \K(\Delta)=\frac{\rkH !}{\rkH_{1}!\cdots \rkH_{k}!}.
\end{equation}
\subsection{}\label{ss;standard2}
Take a pair of weights $\lm,\mu\in\h_{\rkg}^*\subset \t^*_\rkg$
 of $\g_{\rkg}$
such that $\lm-\mu\in P(\Vn)$.
Then there exist integers 
$(\rkH_1,\dots,\rkH_\rkg)\in \Z_{\geq 0}^\rkg$ 
such that $\rkH=\sum_{i=1}^{\rkg}\rkH_i$ and
\begin{equation}\label{eq;partition}
\lm-\mu\equiv\sum_{i=1}^{\rkg}\rkH_i\e_i\quad \text{\rom{mod }}\ 
\C\e,
\end{equation}
where $\e=\e_1+\cdots+\e_n$.
For $\pair\in\h_{\rkg}^*$, we associate an ordered sequence 
of segments
\begin{align}
\Delta_{\pair}:&=
\left(\, [\mu_{1}',\mu_1'
+\rkH_1-1],
\dots,[\mu_{\rkg}',\mu_{\rkg}'+{\rkH}_\rkg-1]\,\right),
\end{align}
where $\mu_i'=\bra \mu+\rho, \e\ch_i\ket$.
We put 
\begin{align}
  \K(\pair)=\K(\Delta_{\lm,\mu}),\quad
  \one_\pair=\one_{\Delta_\pair},\label{eq;one_pair}
\end{align}
where $\one_{\Delta_\pair}$ is as in \eqref{eq;one}.
We call $\K(\lm,\mu)$ a {\it standard module} 
if $\lm+\rho\in P^+_\rkg$.
It is known that
the standard module $\K(\pair)$ has a unique simple quotient,
which is denoted by $\L(\pair)$
(see Theorem \ref{th;USQ}).
%
\subsection{Images of Verma modules}
Our goal in this subsection is the following.
\begin{th}\label{th;Verma_to_standard}
 For $\lm,\mu\in P_{\rkg}$,
there is an isomorphism of $H_{\rkH}$-modules
\begin{equation*}
  \label{eq;verma-standard}
  F_\lm(M(\mu))\cong
\left\{
\begin{array}{ll}
\K(\pair)\quad &\hbox{{ if }}\lm-\mu\in P(\Vn),\\ 
 0                &\hbox{{ otherwise}},
\end{array}\right.
\end{equation*}
where $\K(\pair)$ is given by
 \eqref{eq;one_pair} and \eqref{eq;standard}.
In particular, if $\lm+\rho\in P_\rkg^+$ and $\lm-\mu\in P(\Vn)$,
then
 $F_\lm(M(\mu))$ has a unique simple
quotient. 
\end{th}
To prove Theorem \ref{th;Verma_to_standard},
we prepare some lemmas.
For $\mu\in \h^*_\rkg$,
let $v_\mu$ denote the highest-weight vector of $M(\mu)$.
\begin{lem}\label{lem;asVS}
For $\lm,\mu\in P_{\rkg}$, the natural inclusion 
$(\Vn)_{\lm-\mu}\hookrightarrow (M(\mu)\*\Vn)_{\lm}$
given by $u\mapsto v_\mu\* u$
induces an isomorphism as $W_{\rkH}$-modules
\begin{equation}
  \label{eq;asVS}
 (\Vn)_{\lm-\mu}\buildrel\sim\over\rightarrow F_\lm(M(\mu)).
\end{equation}
In particular  $F_\lm(M(\mu))=0$ unless  $\lm-\mu\in P(\Vn)$.
\end{lem}
\begin{pf}
The lemma follows from the following
fact (known as the tensor product formula):
For any $\mu\in\h^*_{\rkg}$ and any $\g_{\rkg}$-module $Y$
there exists a unique $\g_\rkg$-isomorphism
$$M(\mu)\* Y\buildrel\sim\over\rightarrow
 \hbox{Ind}_{U(\n_+\+\h)}^{U(\g_\rkg)}(\C v_\mu\* Y),$$
which sends $v_\mu\*u$ to $v_\mu\* u$ $(u\in Y)$.
\end{pf}
Recall that
 $\{u_i\}_{i=1,\dots,n}$
is the standard basis of $V_{\rkg}$.
Fix $\lm,\mu\in P_\rkg$ such that $\lm-\mu\in P(\Vn)$
and let $ u_{\pair}\in F_\lm(M(\mu))$ be the image of 
$$
\tilde{u}_{\pair}:
=v_\mu\* u_1^{\* \rkH_1}\* \cdots \* u_{\rkH}^{\* \rkH_{\rkg}}
\in (M(\mu)\* \Vn)_{\lm},$$
where $\rkH_i$ are as in \eqref{eq;partition}.
Let $\zeta_\pair\in\t_\rkH^*$ be the
weight of $\one_\pair$ with respect to the action of $\t_\rkH$:
$$\xi \one_\pair=\bra \zeta_\pair,\xi\ket\one_\pair
\quad (\xi\in\t_\rkH)\quad 
\text{ (see \eqref{eq;one} \eqref{eq;weight_of_one}) }.$$
\begin{lem}\label{lem;weight_vector}
Let $\bar y_i$ denote the image of the operator
$y_i$ $($see \eqref{eq;y_i}$)$ in $\End_\C (F_\lm(M(\mu))$.
Then we have
  $$\bar y_i u_\pair=\bra \zeta_{\lm,\mu}, \e\ch_i\ket u_{\lm,\mu}
\quad (i=1,\dots,\rkH).$$
\end{lem}
{\it Proof.}
Let $\proj: (M(\mu)\*\Vn)_\lm\to F_\lm(M(\mu))$ be the natural surjection.
Then it can be checked that 
$\bar y_i u_\pair
=\proj(v)$,
where 
\begin{align*}
v&=  
\left[\sum_{1 \leq  j<i}\left(r_{ji}+\frac{1}{2\rkg}\right)
     -\sum_{i<j\leq \rkH}  \left(r_{ij}+\frac{1}{2\rkg}\right)
-\pi_i\left(
\frac{1}{2}\sum_{k=1}^{n-1}h_k h^k+
\sum_{\al\in R_n^+} e_{-\al}e_{\al}+\frac{1}{2n}
\right)\right.\\
&\left.
+\frac{1}{2}
\sum_{k=1}^{\rkg-1}\bra \lm+\mu,h^k\ket \pi_i(h_k) -\frac{1}{2}
+\frac{\rkH}{2\rkg}\right] \tilde{u}_{\pair}.
\end{align*}
Here 
$r_{ij}=(\pi_i\*\pi_j)(r)$ ($r$ is as in
\eqref{eq;r-matrix}), and we used
$r_{ij}+r_{ji}=\Omega_{ij}$. 
Now the statement follows from the following formulas
\begin{align*}
&\quad\quad \left( r +\frac{1}{2\rkg} \right) (u_j\* u_k)
=\frac{1}{2}\delta_{jk}(u_j\* u_k)
\quad \text{for }j\leq k,\\
&\left(
\frac{1}{2}\sum_{k=1}^{n-1}h_k h^k+
\sum_{\al\in R_n^+} e_{-\al}e_{\al}+\frac{1}{2n}
\right)u_j
=-\frac{1}{2}(n-2j+1) u_j,\quad
\end{align*}
which are proved by direct calculations.\qed

\medskip
\noindent
{\it Proof of Theorem \ref{th;Verma_to_standard}.}
By Lemma \ref{lem;asVS}, we have that

\smallskip
(i) $u_\pair$ is a cyclic vector of $F_\lm(M(\mu))$,

\smallskip
\noindent
and obviously 

\smallskip
(ii) $wu_\pair=u_\pair$ for all $w\in W_{\rkH_1}\times\cdots\times 
W_{\rkH_{\rkg}}$.

\smallskip\noindent
By  Lemma \ref{lem;weight_vector} and (i)(ii) above,
we have a surjective $H_\rkH$-homomorphism
 $\K(\lm,\mu)\to F_\lm(M(\mu))$ 
which sends $\one_\pair$ to $u_\pair$,
and it is a bijection by Lemma \ref{lem;asVS}. 
\subsection{Images of simple modules}
Next let us suppose that $\rkg=\rkH$ and determine 
the images of simple modules.
\begin{th}\label{th;simple_to_simple}
  Let $\lm+\rho\in P^+_{\rkH}$ and $w\in W_{\rkH}$
be such that $\lm-w\circ\lm\in P(\Vnn)$.
Then we have the following$:$

\medskip
\rom{(i)} If $w$ satisfies
\begin{equation}\label{eq;nonzero_condition}
  \bra w\circ\lm+\rho, h_i\ket\leq 0\quad 
\text{ for any } i\in\{1,\dots,\rkH\} 
\text{ such that } \bra \lm+\rho,h_i\ket=0,
\end{equation}
then
\begin{equation}
   F_\lm(L(w\circ\lm))\cong \L(\lm,w\circ\lm),
\end{equation}
where $\L(\pair)$ is a unique simple quotient of $\K(\pair)$
$($see Theorem \ref{th;USQ}$)$.

\smallskip
\rom{(ii)} 
If $w$ does not satisfy the condition \eqref{eq;nonzero_condition},
then
\begin{equation}
   F_\lm(L(w\circ\lm))=0.
\end{equation}
\end{th}
\begin{rem}\label{rem;equivalent}
For $\lm\in P_\rkH^+-\rho$, let 
$W_{\lm+\rho}$ denote the stabilizer
\begin{equation}
  \label{eq;stabilizer}
  W_{\lm+\rho}=\{w\in W_\rkH\mid w(\lm+\rho)=\lm+\rho\},
\end{equation}
which is a parabolic subgroup of $W_\rkH$.
Let $w_{\scsc L}$ and  $w_{\scsc LR}$
denote the unique
longest element in the coset $W_{\lm+\rho}w$
and $W_{\lm+\rho}w W_{\lm+\rho}$, respectively.
Then the condition \eqref{eq;nonzero_condition}
is equivalent to 
\begin{equation}
  \label{eq;nonzero_condition2}
  w\circ\lm=w_{\scsc L}\circ\lm,
\quad \text{or equivalently}\quad
  w\circ\lm=w_{\scsc LR}\circ\lm.
\end{equation}
Note that $w_{\scsc L}\circ\lm=w_{\scsc LR}\circ\lm$.
\end{rem}
\subsection{}
For $\zeta\in\t^*_{\rkH}$, let 
\begin{equation}
  \label{eq;gamma}
  \gamma_\zeta:Z(H_\rkH)=S(\t_\rkH)^{W_{\rkH}}\to \C
\end{equation}
be the homomorphism given by the evaluation at $\zeta$.
The following is a consequence of Theorem \ref{th;simple_to_simple}
and the classification theorem  of
simple affine Hecke algebra modules
(see Corollary  \ref{cor;parametrization}).
\begin{cor}\label{cor;anysimple}
 Let $\lm+\rho\in P^+_{\rkH}$.
Then any finite-dimensional simple
$H_\rkH$-module with the action of $Z(H_\rkH)$
via $\gamma_{\lm+\rho}$ is isomorphic to
$F_\lm(L(w_{\scsc L}\circ \lm))$ for some $w\in W_\rkH$.
\end{cor}
\begin{rem}
For $c\in\C$, let $t_c$ be the automorphism
of $H_\rkH$ given by
\begin{align*}
  t_c (s_i)&= s_i \ (i=1,\dots,\rkH -1),\quad
t_c(\e_i\ch) =\e_i\ch +c\ (i=1,\dots,\rkH).
\end{align*}
For an $H_\rkH$-module $Y$,
let $Y^c$ denote the $H_\rkH$-module given by the composition
\begin{equation*}
  H_\rkH \overset{t_c}{\to} H_\rkH\to \End_\C (Y).
\end{equation*}
It is known that
any simple $H_\rkH$-module is isomorphic to
\begin{equation*}
{H_\rkH}\*_{ 
 H_{\rkH_1}\* \cdots \* H_{\rkH_k}
} 
 \left( \L(\lm^{(1)},\mu^{(1)})^{c_1}
 \* \cdots\*  \L(\lm^{(k)},\mu^{(k)})^{c_k}\right),
\end{equation*}
for some $\{(\rkH_i,\lm^{(i)}.\mu^{(i)},c_i)\}_{i=1,\dots,k}$.
Here $(\rkH_1,\dots,\rkH_k)\in\Z_{>0}$ is a partition of $\rkH$, 
$c_i$ is a complex number, and 
$(\lm^{(i)},\mu^{(i)})\in P_{\rkH_i}\times P_{\rkH_i}$
satisfying $\lm^{(i)}+\rho\in P_{\rkH_i}^+$
and $\lm^{(i)}-\mu^{(i)}\in P(V_{\rkH_i}^{\*\rkH_i})$
(see \cite{Ch;special bases}).
\end{rem}
\subsection{Proof of Theorem \ref{th;simple_to_simple}}
For $w,y \in W_{\rkg}$ such that $w\leq y$
(where $\leq$ denotes the Bruhat order in $W_n$), 
let $P_{w,y}(q)\in \Z [q]$ denote the Kazhdan-Lusztig polynomial
of the Hecke algebra 
associated to 
$W_{\rkg}$ (see \cite{KL;Coxeter groups,KL;Poincare duality}).
We put for convenience 
$P_{w,y}(q)=0$ for $w\not\leq y$.

The key formula in the following proof of 
Theorem \ref{th;simple_to_simple}
is the following formula
(see Theorem \ref{th;KL1},
Theorem \ref{th;KL2} and Corollary \ref{cor;two_multiplicities}):
\begin{equation}\label{eq;key}
    [\K(\lm,w\circ\lm):\L(\lm,y\circ\lm)]=P_{ w_{LR}, y_{LR} }(1)
=[M(w_{\scsc L}\circ \lm ):L(y_{\scsc L}\circ \lm )].
\end{equation}
Here $\lm\in \h^*_\rkH$ and $w,y\in W_\rkH$ are assumed to satisfy
$\lm+\rho\in P^+_\rkH$ and
$\lm-w\circ\lm,$ $\lm-y\circ\lm\in P(\Vnn)$, and
$[M:N]$ denotes the multiplicity of $N$ in the composition 
series of $M$,
and  $y_{{\scsc LR}}$ (resp. $y_{{\scsc L}}$) denotes the
longest element 
in the coset $W_{\lm+\rho}y W_{\lm+\rho}$ (resp. $yW_{\lm+\rho}$).
 
First we show the following lemma. 
\begin{lem}\label{lem;zero_nonzero}
Let  $\lm+\rho\in P_{\rkH}^+$ and $w\in W_{\rkH}$ 
be such that $\lm-w\circ\lm\in P(\Vnn)$.

\medskip 
\rom{(i)} If $w\in W_{\rkH}$ satisfies the condition 
\eqref{eq;nonzero_condition}, then
$F_\lm(L(w\circ\lm))\neq 0$.

\smallskip
\rom{(ii)} If $w\in W_{\rkH}$ does not satisfy
the condition 
\eqref{eq;nonzero_condition}, then
$F_\lm(L(w\circ\lm))=0$
\end{lem}
\begin{pf}
We first prove (ii).
Suppose that $w$ 
does not satisfy the condition 
\eqref{eq;nonzero_condition}.
Then we can find $s_i\in W_{\lm+\rho}$ such that 
$\bra w\circ\lm+\rho,h_i\ket > 0$.
This inequality means that
 $M(w\circ\lm)$ contains $M(s_i w\circ\lm)$ as a (proper) submodule,
and hence $F_\lm(L(w\circ\lm))$ is a quotient of 
$F_\lm(M(w\circ\lm)/M(s_i w\circ\lm))
\cong \K(\lm,w\circ\lm)/\K(\lm, s_iw\circ\lm)$.
Since $s_i\in W_{\lm+\rho}$, we have
$\lm-s_iw\circ\lm\in P(\Vnn)$ and 
 $\dim \K(\lm,w\circ\lm)=\dim\K(\lm,s_iw\circ\lm)$ by 
\eqref{eq;dim_of_standard}, 
and thus we get $F_\lm(L(w\circ\lm))=0$.

Let us prove (i).
Assume that $w$ satisfies the condition 
\eqref{eq;nonzero_condition}.
Then by Remark \ref{rem;equivalent}, we can assume that
$w$ is the longest element in $W_{\lm+\rho}wW_{\lm+\rho}$.
We can write 
in the Grothendieck group of $\O(\g_\rkg)$ as
\begin{equation}
  \label{eq;Verma_JH}
  M(w\circ\lm)=L(w\circ\lm)
+\sum_{y>w}P_{w,y}(1)
L(y\circ\lm)
\end{equation}
(see  Theorem \ref{th;KL1}), and 
the sum runs over those elements
$y\in W_\rkH$ such that $
y$ is longest in $yW_{\lm+\rho}$ and
$y>w$.
Applying $F_\lm$ to \eqref{eq;Verma_JH}, we have
\begin{equation}
\label{eq;eq}\K(\lm,w\circ\lm)=F_\lm(L(w\circ\lm))
+\sum_{y>w}P_{w, y}(1) F_\lm(L(y\circ\lm))
\end{equation}
in the Grothendieck group of $\Rep(H_\l)$.

Now, let us assume that $F_\lm(L(w\circ\lm))=0$.
Since $[\K(\lm,w\circ\lm):\L(\lm,w\circ\lm)]>0$,
there must be a summand $F_\lm(L(y\circ\lm))$ in the 
right hand side of 
\eqref{eq;eq}
such that
\begin{equation}
  \label{eq;image}
  [F_\lm(L(y\circ\lm)): \L(\lm,w\circ\lm)]> 0.
\end{equation}
Since $F_\lm(L(y\circ\lm))$ is a non-zero quotient of
$F_\lm(M(y\circ\lm))$, we have
\begin{equation}\label{eq;contra}
[\K(\lm,y\circ\lm): \L(\lm,w\circ\lm)]\geq 
 [F_\lm(L(y\circ\mu)): \L(\lm,w\circ\lm)]>0.
\end{equation}
On the other hand, \eqref{eq;key}
implies
$$
[\K(\lm,y\circ\lm): \L(\lm,w\circ\lm)]=P_{y_{LR},w_{LR}}(1)=0
$$
since $l(y_{\scsc LR})\geq l(y)>l(w)=l(w_{\scsc LR})$.
This contradicts \eqref{eq;contra}
and shows that $F_\lm(L(w\circ\lm))\neq 0$.
\end{pf}

Let us complete the  proof of Theorem \ref{th;simple_to_simple}.
Assume that $w\in W_\rkH$ satisfies 
$\lm-w\circ\lm\in P(\Vnn)$ and $w$ is the longest element in 
$W_{\lm+\rho} w W_{\lm+\rho}$
(i.e. $w=w_{\scsc LR}$).
We suppose that $F_\lm(L(w\circ\lm))$ has
a constituent (a simple subquotient)
other than $\L(\lm,w\circ\lm)$,
and will deduce a contradiction.
Such a constituent is isomorphic to
$\L(\lm,y\circ\lm)$ for some $y=y_{\scsc LR}\in W_\rkH$
such that $\lm-y\circ\lm\in P(\Vnn)$
(see Corollary \ref{cor;parametrization}):
\begin{equation}\label{eq;ineq1}
[F_\lm(L(w\circ\lm)):L(\lm,y\circ\lm))]\geq 1.
\end{equation}
Since $F_\lm(L(y\circ\lm))$ is a non-zero quotient of 
$\K(\lm,y\circ\lm)$ by Lemma 
\ref{lem;zero_nonzero},
we have
\begin{equation}\label{eq;ineq2}
[F_\lm(L(y\circ\lm)):L(\lm,y\circ\lm)]\geq 1.
\end{equation}
Combining \eqref{eq;eq}
and inequalities \eqref{eq;ineq1} \eqref{eq;ineq2}, we have
\begin{align*}
&[\K(\lm,w\circ\lm):\L(\lm,y\circ\lm)]\\
&\geq 
[F_\lm(L(w\circ\lm):\L(\lm,y\circ\lm)]+
P_{w,y}(1)[F_\lm(L(y\circ\lm)):L(\lm,y\circ\lm)]\\
&\geq 1+P_{w,y}(1),
\end{align*}
which contradicts \eqref{eq;key} since $y=y_{\scsc LR}$
and $w=w_{\scsc LR}$. $\quad \Box$
\appendix
\section{Some facts from representation theory}
We will review some facts used in this paper.
\subsection{Composition series of Verma modules of $\g_\rkg$}
\label{ss;composition}
For $w,y \in W_{\rkg}$ such that $w\leq y$, 
let $P_{w,y}(q)\in \Z [q]$ denote the Kazhdan-Lusztig polynomial
of the Hecke algebra 
associated to $W_{\rkg}$.
We put for convenience $P_{w,y}(q)=0$ for $w\not\leq y$.
Let $W_{\lm+\rho}:=\{w\in W_\rkg\mid w(\lm+\rho)=\lm+\rho\}$
be the stabilizer.
\begin{th}{\bf{\rom{(\cite{BB,BK}).}}}\label{th;KL1}
Let $\lm+\rho\in P^+_{\rkg}$ and $w\in W_{\rkg}$.
Then
any composition factor of $M(w\circ\lm)$
is isomorphic to $L(y\circ\lm)$ for some $y\in W_\rkg$,
and its multiplicity is given by
\begin{equation}\label{eq;KL2}
   [M(w\circ \lm ):L(y\circ \lm )]=P_{w, y_R}(1)
\end{equation}
where $y_{{\scsc R}}$ is the longest element 
in the right coset $yW_{\lm+\rho}$.
\end{th}
\subsection{Finite-dimensional representations of $H_\rkH$}
We review a classification of finite-dimensional simple
$H_{\rkH}$-modules following \cite{Zel;induced}\ 
in terms of our parameterizations.
(Note that the representation theory of $H_\rkH$ is
related to that of the corresponding affine Hecke algebra by
Lusztig \cite{Lu}.)
\begin{th}\label{th;USQ}
{\rom{(\cite[Theorem 6.1-(a)]{Zel;induced}, see also 
\cite[Theorem 5.2]{Rog})}}
Suppose that $\pair\in\h_{\rkg}^*$ satisfy
 $\lm+\rho\in P^+_\rkg$ and $\lm-\mu\in P(\Vn)$.
Then $\K(\pair)$ has a unique simple quotient,
which is denoted by $\L(\lm,\mu)$.
\end{th}
\begin{lem}\label{lem;isom_type}
{\rom{(\cite[Theorem 6.1-(b)]{Zel;induced})}}
Suppose that $\lm+\rho\in P_\rkg^+$
and $\mu,\eta\in P_\rkg$ satisfy
$\lm-\mu\in P(\Vn)$, $\lm-\eta\in P(\Vn)$.
Then  the following conditions are equivalent:

\medskip
\rom{(i)} $\K(\lm,\mu)\cong\K(\lm,\eta)$.

\smallskip
\rom{(ii)} $\L(\lm,\mu)\cong \L(\lm,\eta)$.

\smallskip
\rom{(iii)} There exists $w\in W_{\lm+\rho}$ such that 
$\eta=w\circ\mu$.
\end{lem}
Let us 
restrict ourselves to the case $\rkg=\rkH$.
In this case, the module
$$\K(\lm,\lm)=H_\rkH\*_{S(\t_\rkH)}\C_{\lm+\rho},$$
is isomorphic to the regular representation as a $\C[W_\rkH]$-module,
where $\C_{\lm+\rho}$ is a one-dimensional $S(\t_\rkH)$-module
determined by the weight $\lm+\rho\in \h^*_\rkH\subset \t^*_\rkH$.
This module is called the principal series representation and
studied e.g. in 
\cite{Ch;unification,Ka,Rog}.
\begin{th}
Let $\lm\in\h^*_\rkH$ be such that
$\lm+\rho\in P^+_{\rkH}\subset \t_{\rkH}^*$.

\medskip
\rom{(i)} \rom{(\cite{Ch;unification,Rog})}
Any finite-dimensional simple module with action of $Z(H_\rkH)$ 
via $\gamma_{\lm+\rho}$ $($where $\gamma_{\lm+\rho}$ is as in
\eqref{eq;gamma}$)$
is a constituent of $\K(\lm,\lm)$.

\smallskip
\rom{(ii)} \rom{(\cite[Theorem 6.1-(c),\ Theorem 7.1]{Zel;induced})}
For $w\in W_\rkH$ such that $\lm-w\circ\lm\in P(\Vnn)$,
the module $\L(\lm,w\circ\lm)$ is a constituent of $\K(\lm,\lm)$,
and
any constituent 
is isomorphic to $\L(\lm,w\circ\lm)$ for some $w\in W_\rkH$
such that $\lm-w\circ\lm\in P(\Vnn)$.
\end{th}
Put 
$$\Para(\lm):=\{w\in W_{\rkH}\mid \lm-w\circ\lm\in P(\Vnn)\}
\subseteq W_{\rkH}$$
and let $\overline{\Para(\lm)}$ denote the image of $\Para(\lm)$ in 
the double coset $W_{\lm+\rho}\backslash W_{\rkH}/W_{\lm+\rho}$.
\begin{cor}\label{cor;parametrization}
 Let $\lm+\rho\in P^+_{\rkH}$.
Then there exists a one to one
correspondence between
$\overline{\Para(\lm)}$ and
the set of  equivalent classes of finite-dimensional simple 
$H_\rkH$-modules with the action of $Z(H_\rkH)$
via $\gamma_{\lm+\rho}$, which is
given by
$W_{\lm+\rho}wW_{\lm+\rho}\mapsto
\L(\lm,w\circ\lm).$
\end{cor}
\subsection{Multiplicity formulas}
Let us recall the multiplicity formula for 
(degenerate) affine Hecke algebras of $GL_\rkH$.
Zelevinsky conjectured in \cite{Zel;p-adic KL} that
the multiplicity of simple modules in the composition series of
standard modules is given in terms of
the intersection cohomologies concerning the quiver variety.
He proved in  \cite{Zel;two remarks} 
that these intersection cohomologies are
expressed by the Kazhdan-Lusztig polynomials of
the symmetric group.
Zelevinsky's conjecture was proved by 
Ginzburg\cite{Gi} (see also \cite{CG}) in more general situations.
The result is rephrased as follows:
\begin{th}[\cite{Gi}]\label{th;KL2}
  Suppose that $\lm+\rho\in P^+_{\rkH}$ 
and $y,w\in W_{\rkH}$ satisfy
$\lm-w\circ\lm\in P(\Vnn)$ and  $\lm-y\circ\lm\in P(\Vnn)$. Then
\begin{equation}
    [\K(\lm,w\circ\lm):\L(\lm,y\circ\lm)]=P_{ w_{LR}, y_{LR} }(1),
  \end{equation}
where $w_{{\scsc LR}}$ and $y_{{\scsc LR}}$ denote the longest elements
in the double coset $W_{\lm+\rho} w W_{\lm+\rho}$ and 
$W_{\lm+\rho} y W_{\lm+\rho}$
respectively.
\end{th}
Combining Theorem \ref{th;KL1} and Theorem \ref{th;KL2}
we have the following identity.
\begin{cor}\label{cor;two_multiplicities}
Assume that $\lm+\rho\in P^+_{\rkH}$ and that $w,y\in W_{\rkH}$ satisfy
$\lm-w\circ\lm\in P(\Vnn)$ and $\lm-y\circ\lm\in P(\Vnn)$.
Then we have
\begin{equation}\label{eq;two_multiplicities}
[\K(\lm,w\circ\lm):\L(\lm,y\circ\lm)]
=[M(w_L\circ\lm):L(y_L\circ\lm)],
\end{equation}
where $w_{{\scsc L}}$ and $y_{{\scsc L}}$ denote the longest elements
in the left coset $W_{\lm+\rho} w$ and $W_{\lm+\rho} y$
respectively.
\end{cor}

\end{document}